\documentclass[12pt]{article}

\usepackage{amsfonts}
\usepackage{amssymb}
\usepackage{amsfonts,amssymb, amsmath, latexsym}
\usepackage{pst-all}
\usepackage{CJK}
\usepackage{mathrsfs}

\usepackage{amsfonts,amssymb, amsmath, latexsym}
\usepackage{pst-all}
\usepackage{CJK}
\usepackage[left,pagewise,mathlines]{lineno}
\usepackage{cite}

\newtheorem{thm}{Theorem}

\newtheorem{lem}[thm]{Lemma}

\newtheorem{claim}{Claim}

\newenvironment{pf}[1][Proof]{\noindent\textbf{#1.} }{\hfill\rule{1mm}{2mm}}

\def\g{\gamma}

\begin{document}

\title{The $2$-Domination and $2$-Bondage Numbers\\ of Grid Graphs\thanks{The work was
supported by NNSF of China (No.10711233) and the Fundamental Research Fund of NPU (No. JC201150)}}

\author
{ You Lu\\
{\small Department of Applied Mathematics}\\
{\small Northwestern Polytechnical University}\\
{\small Xi'an Shaanxi 710072, P. R. China}\\
{\small  Email: luyou@nwpu.edu.cn}\\ \\
 Jun-Ming Xu\footnote{Corresponding author: xujm@ustc.edu.cn} \\
{\small School of Mathematical Science}\\
{\small University of Science and Technology of China}   \\
{\small Wentsun Wu Key Laboratory of CAS}\\
{\small Hefei, Anhui, 230026, P. R. China}\\
{\small Email: xujm@ustc.edu.cn}  \\
}
\date{}
\maketitle

\begin{abstract}
Let $p$ be a positive integer and $G=(V,E)$ be a simple graph. A
subset $D\subseteq V$ is a $p$-dominating set if each vertex not in
$D$ has at least $p$ neighbors in $D$. The $p$-domination number
$\g_p(G)$ is the minimum cardinality  among all $p$-dominating sets
of $G$. The $p$-bondage number $b_p(G)$ is the cardinality of a
smallest set of edges whose removal from $G$ results in a graph with
a $p$-domination number greater than the $p$-domination number of
$G$. In this note we determine the $2$-domination number $\g_2$ and
$2$-bondage number $b_2$ for the grid graphs $G_{m,n}=P_m\times P_n$
for $2\leq m\leq 4$.
\\



 \end{abstract}

\section{Induction}

Let $G =(V(G),E(G))$ be an undirected simple graph with vertex set
$V(G)$ and edge set $E(G)$. The {\it open neighborhood} and the {\it
degree} of a vertex $v\in V(G)$ are denoted by $N_G(v)=\{u\ |\ uv\in
E(G)\}$ and $deg_G(v)=|N_G(v)|$, respectively. For a subset
$S\subset V(G)$, the subgraph induced by $V(G)\setminus S$ is
denoted by $G - S$. For any $B\subseteq E(G)$, we use $G-B$ to
denote the subgraph with vertex set $V(G)$ and edge set
$E(G)\setminus B$. For convenience, for $x\in V(G)$ and $uv\in
E(G)$, we denote $G-\{x\}$ and $G-\{uv\}$ by $G-x$ and $G-uv$,
respectively. In this paper, we follow \cite{hs981,hs982} for
graph-theoretical terminology and notation not defined here.

Let $p$ be a positive integer. In \cite{fj851,fj852}, Fink and
Jacobson introduced the concept of $p$-domination. A subset $D$ of
$V(G)$ is a {\it $p$-dominating set} of $G$ if for every vertex
$v\in V(G)$, $|D\cap N_G(v)|\geq p$. The {\it $p$-domination number}
$\g_p(G)$ is the minimum cardinality among the $p$-dominating sets
of $G$. Any $p$-dominating set of $G$ with cardinality $\g_p(G)$
will be called a $\g_p(G)$-set.  For any $S,T\subseteq V(G)$, $S$
$p$-dominates $T$ in $G$ if for every vertex $v\in T$, $|S\cap
N_G(v)|\geq p$. Notice that the $1$-dominating set is a classical
dominating set, and so $\g_1(G)=\g(G)$. The $p$-domination number
has received much research attention (see, for example,
\cite{bcf05,bcv06,c06,cr90,dghpv11,f85}).

In particular, when $p=2$, Fink and Jacobson~\cite{fj851}
established $\gamma_2(T)\geq(n+1)/2$ for a tree $T$ of order $n$ and
Blidia et al~\cite{bcv06b} showed $\gamma_2(G)\leq(n+\gamma)/2$ for
a graph of order $n$, minimum degree $\delta\geq2$, and domination
number $\gamma$. Since then, the study on $2$-domination number has
received much attention (see, for example, \cite{bcv06b,c06,
dlpw10,hv07,hv08a,hv08b,k10,s09,v08}).

As a measurement of the stability of $p$-domination in a graph under
edge removal, Lu and Xu \cite{lx10} introduced the $p$-bondage
number $b_p(G)$. The {\it $p$-bondage number} of $G$ is the minimum
cardinality among all edge subsets $B \subseteq E(G)$ such that
$\g_p(G-B) > \g_p(G)$. The case $p=1$ leads to the usual bondage
number $b(G)$, which is introduced by Fink et al.~\cite{fjkr90} and
further study for example in
\cite{cd06,frv03,hr92,hl94,hx06,ksk05,t97}.

The notation $[i]$ denotes the set $\{1,2,\cdots,i\}$. The Cartesian
product $G_{m,n}=P_m\times P_n$ of two paths is the {\it grid graph}
with vertex set $V(G_{m,n})=[m]\times [n]$, where two vertices
$(i,j)$ and $(i',j')$ are adjacent if and only if $|i-i'|+|j-j'|=1$.
Notice that there are many research articles on the $\g_p(G_{m,n})$
for $p=1$. In 1983, Jacobson and Kinch \cite{jk83} established the
exact values of $\g(G_{m,n})$ for $2\leq m\leq 4$ which are the
first results on the domination number of grids. In 1993, Chang and
Clark~\cite{cc93} found those of $\g(G_{m,n})$ for $m=5$ and $6$.
Fischer found those of $\g(G_{m,n})$ for $m\leq 21$ (see
Gon\c{c}alves et al.\cite{gprt11}).
Recently, Gon\c{c}alves et al.\cite{gprt11} finished the computation
of $\g(G_{m,n})$ when $24\leq m\leq n$. In \cite{hcx09}, the authors
determined $b_1(G_{m,n})$ for $2\leq m\leq 4$.
Until now, however, no research has been done on calculating the
values of $\g_p(G_{m,n})$ and $b_p(G_{m,n})$ for $p\geq 2$. In this
paper, we will determine the values of $\g_2(G_{m,n})$ and
$b_2(G_{m,n})$ for $2\leq m\leq 4$.

The rest of this paper is organized as follows. In Section 2, we
give some lemmas and determine the $2$-domination and $2$-bondage
numbers of $G_{2,n}$. We determine the $2$-domination and
$2$-bondage numbers of $G_{3,n}$ and $G_{4,n}$ in Section 3 and
Section 4, respectively.

\section{Preliminaries}
\begin{lem}\label{l1}
Every $2$-dominating set of a graph contains all vertices of degree
one.
\end{lem}

Throughout this paper, we will use the following notations. Let $D$
be a $2$-dominating set $D$ of $G_{m,n}$. For any $j\in [n]$, let
 $$
 C_j=\{(i,j)\ |\ i\in [m]\},\ \ \ \mathscr{C}_j(D)=D\cap C_j\ \ \ \mbox{and}\ \ \ c_j(D)=|\mathscr{C}_j(D)|.
 $$
The sequences
$(\mathscr{C}_1(D),\mathscr{C}_2(D),\cdots,\mathscr{C}_n(D))$ and
$(c_1(D),c_2(D),\cdots,c_n(D))$ will be called $2$-dominating set
and $2$-domination number sequences of $G_{n,m}$, respectively.
Since $0\leq c_j(D)\leq m$ for each $j\in [n]$, let
 $$
 N_i(D)=|\{c_j(D)\ |\ c_j(D)=i,\ \  j\in [n]\}|,\ \ \ \mbox{where $i=0,1,\cdots,m$}.
 $$

Note that every vertex in $\mathscr{C}_2(D)$ dominates exactly one
vertex of $C_1$, and so the column $C_1$ is $1$-dominated by
$\mathscr{C}_1(D)$. With the symmetry of $C_1$ and $C_n$, we have
\begin{lem}\label{l2}
For any $m,n\geq 1$, $c_1(D)\geq \lceil\frac{m}{3}\rceil$ and
$c_n(D)\geq \lceil\frac{m}{3}\rceil$.
\end{lem}

\begin{lem}\label{l3}
Let $m,n\geq 3$ and $2\leq j\leq n-1$. If
$\mathscr{C}_j(D)=\{(i,j)\}$ with $i\in \{1,2,m-1,m\}$, then
$c_{j-1}(D)+c_{j+1}(D)\geq 2m-3$.
\end{lem}

\begin{pf} By the symmetry of $G_{m,n}$, we only need to prove the cases
$i=1$ and $2$. Note that $C_j$ is $2$-dominated by
$\mathscr{C}_{j-1}(D)\cup \mathscr{C}_{j}(D)\cup
\mathscr{C}_{j+1}(D)$ in $G_{m,n}$. So $C_{j-1}\setminus
\{(1,j-1),(2,j-1), (3,j-1)\}\subseteq D$ and $C_{j+1}\setminus
\{(1,j+1),(2,j+1),(3,j+1)\}\subseteq D$.

If $i=1$, then $\{(3,j-1),(3,j+1)\}\subseteq D$ and $|D\cap
\{(2,j-1),(2,j+1)\}|\geq 1$. Hence
 $$
 c_{j-1}(D)+c_{j+1}(D)
                \geq (m-3)+(m-3)+2+1=2m-3.
 $$

If $i=2$, to $2$-dominate $(3,j)$, $|D\cap \{(3,j-1),(3,j+1)\}|\geq
1$. Since $(1,j)\notin D$, to $2$-dominate $(1,j-1)$ and $(1,j+1)$,
we have $|D\cap \{(1,j-1),(2,j-1)\}|\geq 1$ and $|D\cap
\{(1,j+1),(2,j+1)\}|\geq 1$. Thus
 $$
 c_{j-1}(D)+c_{j+1}(D)\geq (m-3)+(m-3)+1+1+1=2m-3.
 $$
The lemma follows.\end{pf}

For any $2$-dominating set $D$ of $G_{m,n}$, if there is some $j\in
[n]$ with $c_j(D)=0$, then $c_{j-1}(D)=c_{j+1}(D)=m$. This forces
        \begin{equation}\label{e2.1}
              N_0(D) \leq \left\{
                    \begin{array}{l}
                                  0 \hspace{2.15cm}\ {\rm if}\ \ N_m(D)=0; \\
                                  N_m(D)-1 \hspace{0.5cm} {\rm if}\ \ N_m(D)\geq 1.
                    \end{array}
           \right.
        \end{equation}

\begin{lem}\label{l4}
If $G_{m,n}$ contains a $\g_2(G_{m,n})$-set $S$ with $N_0(S)\neq 0$,
then there must be another $\g_2(G_{m,n})$-set $D$ such that
$N_0(D)=0$ and $N_m(D)\neq 0$.
\end{lem}

\begin{pf} From (\ref{e2.1}) and $N_0(S)\neq 0$, we know that
$N_m(S)-N_0(S)\geq 1$. Denote $\min\{j\ |\ c_j(S)=0, j\in [n]\}$ by
$t$. To $2$-dominate $C_t$, $c_{t-1}=m$ and $c_{t+1}=m$, and so
$t\geq 2$. By the choice of $t$, we must have $\mathscr{C}_{t-2}\neq
\emptyset$ (if $t\geq 3$), and let $(i,t-2)\in \mathscr{C}_{t-2}$.
Clearly, subset $S_1=(S\setminus \{(i,t-1)\})\cup \{(i,t)\}$ is a
$\g_2(G_{m,n})$-set with $N_0(S_1)=N_0(S)-1$ and
$N_m(S_1)=N_m(S)-1$. By using recursively the above operation
$N_0(S)-1$ times,  we can obtain  a $\g_2(G_{m,n})$-set $D$  from
$S_1$, and $D$ satisfies that $N_0(D)=N_0(S)-N_0(S)=0$ and
$N_m(D)=N_m(S)-N_0(S)\geq 1$. The result holds. \end{pf}

Applying the above lemmas, the values of $\g_2(G_{2,n})$ and
$b_2(G_{2,n})$ can be easily obtained.

\begin{thm}\label{t1}
For any positive integer $n\geq 2$, $\g_2(G_{2,n})=n$ and $b_2(G_{2,n})=1$.
\end{thm}

\begin{pf}
We first prove $\g_2(G_{2,n})=n$. Note that vertex subset
 $$ S=\left\{
                    \begin{array}{l}
                                  \cup_{i=1}^k\{(2i-1,1),(2i,2)\} \hspace{2.2cm} {\rm if}\  n=2k \\
                                  \cup_{i=1}^k\{(2i-1,1),(2i,2)\}\cup \{(n,1)\} \hspace{0.35cm} {\rm if}\  n=2k+1
                    \end{array}
           \right.$$
 is a $2$-dominating set of $G_{2,n}$. So $\g_2(G_{2,n})\leq |S|=n$.

Let $D$ be a $\g_2(G_{2,n})$-set, then
$N_1(D)+2N_2(D)=|D|=\g_2(G_{2,n})\leq n$ and
$N_0(D)+N_1(D)+N_2(D)=n$. This forces $N_2(D)\leq N_0(D)$. By
(\ref{e2.1}), $N_0(D)=N_2(D)=0$, and so the $2$-domination number
sequence $$(c_1(D),c_2(D),\cdots,c_n(D))=(1,1,\cdots,1).$$ Hence
$\g_2(G_{2,n})=N_1(D)=n$. By the arbitrariness of $D$, all
domination number sequences of $G_{2,n}$ is $(1,1,\cdots,1)$.

Next we prove $b_2(G)=1$. Let $H=G_{2,n}-(1,1)(2,1)$ and $D'$ be a
$\g_2(H)$-set. By Lemma \ref{l1}, $(1,1)\in D'$ and $(2,1)\in D'$
since $(1,1)$ and $(2,1)$ are two vertices with degree one in $H$.
Since $D'$ is a $2$-dominating set of $G_{2,n}$ with $c_1(D')=2$, we
know that $D'$ isn't a $\g_2(G_{2,n})$-set, and so
  $
\g_2(H)=|D'|> \g_2(G_{2,n}),
  $
which implies $b_2(G_{2,n})=1$.\end{pf}

\section{The grid graphs $G_{3,n}$}
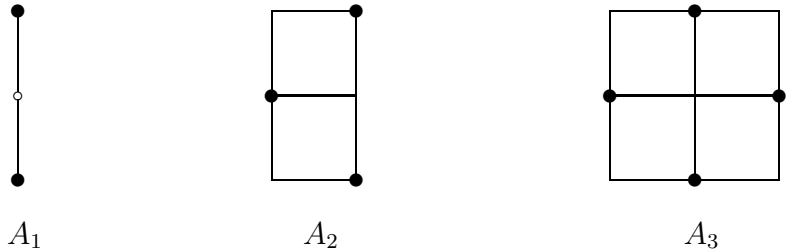
\begin{figure}[h]
\begin{center}
\psset{unit=0.8pt}

\begin{pspicture}(400,200)(0,50)

\put(300,90){\line(1,0){80}}\put(300,170){\line(1,0){80}}
\put(300,90){\line(0,1){80}}\put(380,90){\line(0,1){80}}
\put(300,130){\line(1,0){80}}\put(340,90){\line(0,1){80}}
\put(300,130){\circle*{6}}\put(340,170){\circle*{6}}
\put(340,90){\circle*{6}}\put(380,130){\circle*{6}}
\put(335,60){$A_3$}

\put(20,90){\line(0,1){38}}\put(20,132){\line(0,1){38}}
\put(20,130){\circle{4}}
\put(20,90){\circle*{6}}\put(20,170){\circle*{6}}
\put(15,60){$A_1$}

\put(140,90){\line(1,0){40}}\put(140,170){\line(1,0){40}}
\put(140,90){\line(0,1){80}}\put(180,90){\line(0,1){80}}
\put(140,130){\line(1,0){40}}
\put(140,130){\circle*{6}}\put(180,170){\circle*{6}}
\put(180,90){\circle*{6}}
\put(155,60){$A_2$}
 \end{pspicture}
\caption{\label{f1} \footnotesize Blocks $A_1$, $A_2$ and $A_3$ for constructing $2$-dominating sets of $G_{3,n}$}
\end{center}
\end{figure}

In this section, we will propose the values of $\g_2(G_{3,n})$ and
$b_2(G_{3,n})$. We will construct $2$-dominating set of $G_{3,n}$ by
concatenating the blocks $A_1$, $A_2$ and $A_3$ of Figure \ref{f1},
where the concept of concatenation was introduced by Chang et al. in
\cite{cc93,cch94}.

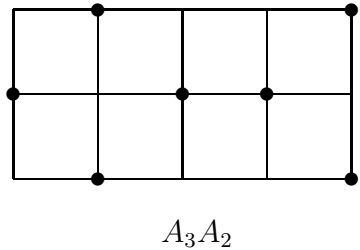
\begin{figure}[h]
\begin{center}
\psset{unit=0.8pt}

\begin{pspicture}(400,200)(0,50)

\put(120,90){\line(1,0){160}}\put(120,130){\line(1,0){160}}
\put(120,170){\line(1,0){160}}
\put(120,90){\line(0,1){80}}\put(160,90){\line(0,1){80}}
\put(200,90){\line(0,1){80}}\put(240,90){\line(0,1){80}}
\put(280,90){\line(0,1){80}}
\put(120,130){\circle*{6}}\put(160,170){\circle*{6}}
\put(160,90){\circle*{6}}\put(200,130){\circle*{6}}
\put(240,130){\circle*{6}}\put(280,90){\circle*{6}}
\put(280,170){\circle*{6}}
\put(190,60){$A_3A_2$}

 \end{pspicture}
\caption{\label{f2} \footnotesize Block $A_3A_2$ by concatenating $A_3$ and $A_2$}
\end{center}
\end{figure}

We explain the meaning of {\it concatenation} by an example: if we
concatenate the $3\times 3$ block $A_3$ and the  $3\times 2$ block
$A_2$ then we obtained the $3\times 5$ block $A_3A_2$ of Figure
\ref{f2}.  For $t\geq 0$, we use
$(A_3A_2)^t=(A_3A_2)(A_3A_2)\cdots(A_3A_2)$ to denote the
concatenation of $A_3A_2$ with itself $t$ times.

Clearly,  the $r+1$ black vertices in $A_r$ constitute a
$2$-dominating set of $G_{3,r}$ for each $r\in \{1,2,3\}$. So the
set of black vertices in $A_3A_2$ is  a $2$-dominating set of
$G_{3,5}$ we wish to construct.

\begin{thm}\label{t2}
For any positive integer $n\geq 1$,
 $
 \g_2(G_{3,n})=\lceil\frac{4n}{3}\rceil.
 $
\end{thm}

\begin{pf} If $n=1,2$ or $3$, then it is trivial
by the definition of $\g_2(G)$-set. In the following, assume that
$n\geq 4$.

Let $n=3k+r$, where $1\leq r\leq 3$. Clearly, the set of black
vertices in $(A_3)^kA_r$ is a $2$-dominating set of $G_{3,n}$, and
so $\g_2(G_{3,n})\leq 4k+(r+1)=\lceil\frac{4n}{3}\rceil$.

Let $D$ be a $\g_2(G_{3,n})$-set. We consider two cases.

\noindent\textbf{Case 1}\ \ $N_0(D)=0$. Then we have
\begin{equation}\label{e3.2}
    N_1(D)+N_2(D)+N_3(D)=n
\end{equation}
and
\begin{equation}\label{e3.3}
    N_1(D)+2N_2(D)+3N_3(D)=|D|=\g_2(G_{3,n})\leq \lceil\frac{4n}{3}\rceil,
\end{equation}
and so $N_1(D)\geq 1$. Let $j\in [n]$ be an integer with $c_j(D)=1$.
If $j=1$, to $2$-dominate $C_1$, we have $c_2(D)\geq 2$; If $j=n$,
then $c_{n-1}(D)\geq 2$ by the symmetry of $C_1$ and $C_n$; If
$2\leq j\leq n-1$, by Lemma \ref{l3}, at least one of
$c_{j-1}(D)\geq 2$ and $c_{j+1}(D)\geq 2$ is true, that is, there
are at most two 1's between two adjacent $a$'s and $b$'s in the
$2$-domination number sequence $(c_1(D),c_2(D),\cdots,c_n(D))$,
where $a\in \{2,3\}$ and $b\in \{2,3\}$. Hence
\begin{equation}\label{e3.4}
    N_1(D)\leq 1+2[(N_2(D)+N_3(D))-1]+1=2(N_2(D)+N_3(D)).
\end{equation}
Combined with (\ref{e3.2}), (\ref{e3.3}) and (\ref{e3.4}), we must have
 $$ \begin{array}{rl}
 \lceil\frac{4n}{3}\rceil\geq\g_2(G_{3,n})&=N_1(D)+2N_2(D)+3N_3(D)\\
                                       &=n+\frac{1}{3}[2(N_2(D)+N_3(D))+N_2(D)+N_3(D)]+N_3(D)\\
                                       &\geq n+\frac{1}{3}(N_1(D)+N_2(D)+N_3(D))+N_3(D)\\
                                       &=\frac{4n}{3}+N_3(D),
 \end{array}$$
which implies $\g_2(G_{3,n})=\lceil\frac{4n}{3}\rceil$ and $N_3(D)=0$.

\noindent\textbf{Case 2}\ \ $N_0(D)\neq 0$. By Lemma \ref{l4},
$G_{3,n}$ contains a $\g_2(G_{3,n})$-set $S$ with $N_0(S)=0$ and
$N_3(S)\neq 0$. By $N_0(S)=0$ and \textbf{Case 1}, we must have
$N_3(S)=0$, which contradicts with $N_3(S)\neq 0$. This completes
the proof.\end{pf}

From the proof of Theorem \ref{t2}, we can know that for any
$\g_2(G_{3,n})$-set $D$, $N_0(D)=N_3(D)=0$, that is, $c_j(D)=1$ or
$2$ for each $j\in [n]$.

\begin{thm}\label{t3}
For any positive integer $n\geq 2$,
 $$
 b_2(G_{3,n})=\left\{\begin{array}{ll}
                                                                   2       & \mbox{{\rm if} {\rm $n\equiv 1$ (mod $3$)}} \\
                                                                   1       & {\rm otherwise.}\\
                      \end{array}
               \right.
 $$
\end{thm}

\begin{pf} It is trivial for $n=2$ or $3$. In the following, assume that $n\geq 4$, and let $n=3k+r$, $1\leq r\leq 3$.
We consider two cases.

\noindent\textbf{Case 1}\ \  $r=3$ or $2$.

Let $G=G_{3,n}-e$ and $S$ be a $\g_2(G)$-set, where $e=(1,1)(1,2)$.
Clearly, $S$ is a $2$-domination of $G_{3,n}$. In the following we
will show $S$ isn't a $\g_2(G_{3,n})$-set, which implies
$b_2(G_{3,n})=1$.

Supposed that $S$ is a $\g_2(G_{3,n})$-set, then we must have
 $$
 N_2(S)=\lceil\frac{n}{3}\rceil=k+1
 $$
  and
 $$
 N_1(S)=|S|-2N_2(S)=\lceil\frac{4n}{3}\rceil-2(k+1)=(4k+r+1)-2k-2=2k+r-1.
 $$
By Lemma \ref{l1}, we have $(1,1)\in S$ and $|\{(2,1),(3,1)\}\cap
S|\geq 1$ since $deg_G((1,1))=1$ and $deg_G((3,1))=2$, and so
$c_1(S)=|C_1\cap S|\geq 2$. By Lemma \ref{l3} and $N_2(S)=k+1$,
$N_1(S)$ is at most $
 2(N_2(S)-1)+1\ (=2k+1)
 $
in the domination number sequence $(c_1(S),c_2(S),\cdots,c_n(S))$.

When $r=3$,  $N_1(S)=2k+2> 2(N_2(S)-1)+1$ contradicts with $N_1(S)\leq 2(N_2(S)-1)+1$.

When $r=2$,  $N_1(S)=2(N_2(S)-1)+1=2k+1$, and so
$(c_1(S),c_2(S),c_3(S))=(2,1,1)$. Since $deg_G((1,2))=2$,  we have
either $(1,2)\in S$ or $\{(2,2),(1,3)\}\subseteq S$. If $(1,2)\in
S$, then $\mathscr{C}_2=\{(1,2)\}$. To $2$-dominate $(3,2)$,
$(3,1)\in S$ and $(3,3)\in S$. Thus $\mathscr{C}_1=\{(1,1),(3,1)\}$,
$\mathscr{C}_2=\{(1,2)\}$ and $\mathscr{C}_3=\{(3,3)\}$. Clearly,
$S$ cannot $2$-dominate vertex $(2,2)$ in $G_{3,n}$, a
contradiction. If $\{(2,2),(1,3)\}\subseteq S$, then
$\mathscr{C}_2=\{(2,2)\}$ and $\mathscr{C}_3=\{(1,3)\}$. It is
obvious that vertex $(3,3)$ cannot be $2$-dominated by $S$ in
$G_{3,n}$, a contradiction.

\noindent\textbf{Case 2}\ \ $r=1$.

We first show that $b_2(G_{3,n})\geq 2$, that is
$\g_2(G_{3,n}-e)=\g_2(G_{3,n})$ for any $e=(i,j)(i',j')\in
E(G_{3,n})$. By the symmetry of $G_{3,n}$, assume $1\leq i \leq 2$
and $1\leq j\leq \lfloor\frac{n}{2}\rfloor$. Let
$t=\lfloor\frac{j}{3}\rfloor$. Note that the set of black vertices,
denoted by $D'$, in block $(A_3)^tA_1(A_3)^{k-t}$ is a
$\g_2(G_{3,n})$-set satisfying one end of edge $e$ is in $D'$ and
another end is dominated by at least three vertices of $D'$. That is
to say, $D'$ is also $2$-dominating set of $G_{3,n}-e$. So we can
obtain $\g_2(G_{3,n}-e)= \g_2(G_{3,n})$ from that $G_{3,n}-e$ is a
spanning subgraph of $G_{3,n}$.

To the end, we merely prove $b_2(G_{3,n})\leq 2$. Let
$e_1=(1,1)(2,1)$, $e_2=(2,1)(3,1)$ and $H=G_{3,n}-\{e_1,e_2\}$. Let
$S$ be a $\g_2(H)$-set, then $S$ is a $2$-dominating set of
$G_{3,n}$. By Lemma \ref{l1} and
$deg_H((1,1))=deg_H((2,1))=deg_H((3,1))=1$, $S$ contains vertices
$(1,1)$, $(2,1)$ and $(3,1)$. Thus $c_1(S)=3$, and so $S$ isn't a
$\g_2(G_{3,n})$-set. This forces $\g_2(H)=|S|>\g_2(G_{3,n})$, which
implies $b_2(G_{3,n})\leq 2$. The proof is completed.\end{pf}

\section{The grid graphs $G_{4,n}$}

In this section, we will present the values of the $2$-domination
and $2$-bondage numbers of $G_{4,n}$. Throughout this section $n\geq
3$.

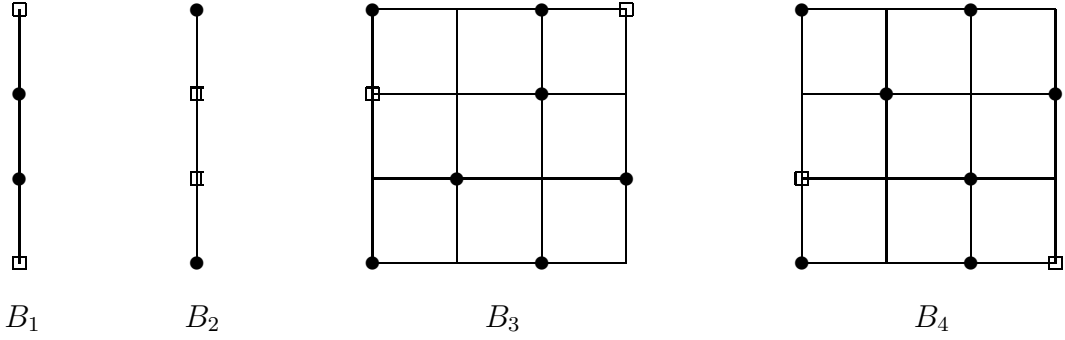
\begin{figure}[h]
\begin{center}
\psset{unit=0.8pt}

\begin{pspicture}(430,200)(0,50)

\put(-30,90){\line(0,1){40}}\put(-30,130){\line(0,1){40}}\put(-30,170){\line(0,1){40}}
\put(-30,170){\circle*{6}}\put(-30,130){\circle*{6}}

\put(-33,87){\line(1,0){6}}\put(-33,93){\line(1,0){6}}
\put(-33,87){\line(0,1){6}}\put(-27,87){\line(0,1){6}}

\put(-33,207){\line(1,0){6}}\put(-33,213){\line(1,0){6}}
\put(-33,207){\line(0,1){6}}\put(-27,207){\line(0,1){6}}

\put(-37,60){$B_1$}

\put(54,90){\line(0,1){40}}\put(54,130){\line(0,1){40}}\put(54,170){\line(0,1){40}}
\put(54,90){\circle*{6}}\put(54,210){\circle*{6}}

\put(51,127){\line(1,0){6}}\put(51,133){\line(1,0){6}}
\put(51,127){\line(0,1){6}}\put(56,127){\line(0,1){6}}

\put(51,167){\line(1,0){6}}\put(51,173){\line(1,0){6}}
\put(51,167){\line(0,1){6}}\put(56,167){\line(0,1){6}}

\put(48,60){$B_2$}

\put(137,90){\line(1,0){120}}\put(137,130){\line(1,0){120}}\put(137,170){\line(1,0){120}}\put(137,210){\line(1,0){120}}
\put(137,90){\line(0,1){120}}\put(177,90){\line(0,1){120}}\put(217,90){\line(0,1){120}}\put(257,90){\line(0,1){120}}

\put(137,90){\circle*{6}}\put(137,210){\circle*{6}}\put(177,130){\circle*{6}}
\put(217,90){\circle*{6}}\put(217,170){\circle*{6}}\put(217,210){\circle*{6}}
\put(257,130){\circle*{6}}

\put(254,207){\line(1,0){6}}\put(254,213){\line(1,0){6}}
\put(254,207){\line(0,1){6}}\put(260,207){\line(0,1){6}}

\put(134,167){\line(1,0){6}}\put(134,173){\line(1,0){6}}
\put(134,167){\line(0,1){6}}\put(140,167){\line(0,1){6}}

\put(190,60){$B_3$}

\put(340,90){\line(1,0){120}}\put(340,130){\line(1,0){120}}\put(340,170){\line(1,0){120}}\put(340,210){\line(1,0){120}}
\put(340,90){\line(0,1){120}}\put(380,90){\line(0,1){120}}\put(420,90){\line(0,1){120}}\put(460,90){\line(0,1){120}}

\put(340,90){\circle*{6}}\put(340,210){\circle*{6}}\put(380,170){\circle*{6}}
\put(420,90){\circle*{6}}\put(420,130){\circle*{6}}\put(420,210){\circle*{6}}
\put(460,170){\circle*{6}}

\put(457,87){\line(1,0){6}}\put(457,93){\line(1,0){6}}
\put(457,87){\line(0,1){6}}\put(463,87){\line(0,1){6}}

\put(337,127){\line(1,0){6}}\put(337,133){\line(1,0){6}}
\put(337,127){\line(0,1){6}}\put(343,127){\line(0,1){6}}

\put(393,60){$B_4$}

 \end{pspicture}
\caption{\label{f3} \footnotesize Blocks $B_1$, $B_2$, $B_3$ and $B_4$ for constructing $2$-dominating set of $G_{4,n}$}
\end{center}
\end{figure}

To construct the $2$-dominating set of $G_{4,n}$, we need the
following four blocks of Figure \ref{f3}.

For each $i\in [4]$, it is obvious that  the squared vertices cannot
be $2$-dominated but can be $1$-dominated by the black vertices in
$B_i$ of Figure \ref{f3}. Thus it is easy to see that the set of
black vertices, denoted by $D_n$, in the block $T_n$ is a
$2$-dominating set of $G_{4,n}$, where the block
 $$
 T_n=            \left\{\begin{array}{ll}
                             B_1(B_3B_4)^{\frac{k}{2}}B_2B_1                       & \mbox{\hspace{0.2cm} if $n=4k+3$ and $k$ is even} \\
                             B_1(B_3B_4)^{\frac{k}{2}}B_2B_1B_2                    & \mbox{\hspace{0.2cm} if $n=4k+4$ and $k$ is even} \\
                             B_1(B_3B_4)^{\frac{k}{2}}B_2B_1B_2B_1                 & \mbox{\hspace{0.2cm} if $n=4k+5$ and $k$ is even} \\
                             \vspace{10pt}
                             B_1(B_3B_4)^{\frac{k}{2}}B_2B_1B_2B_1B_2              & \mbox{\hspace{0.2cm} if $n=4k+6$ and $k$ is even} \\
                             B_1(B_3B_4)^{\frac{k-1}{2}}B_3B_2B_1                  & \mbox{\hspace{0.2cm} if $n=4k+3$ and $k$ is odd} \\
                             B_1(B_3B_4)^{\frac{k-1}{2}}B_3B_2B_1B_2               & \mbox{\hspace{0.2cm} if $n=4k+4$ and $k$ is odd} \\
                             B_1(B_3B_4)^{\frac{k-1}{2}}B_3B_2B_1B_2B_1            & \mbox{\hspace{0.2cm} if $n=4k+5$ and $k$ is odd} \\
                             B_1(B_3B_4)^{\frac{k-1}{2}}B_3B_2B_1B_2B_1B_2         & \mbox{\hspace{0.2cm} if $n=4k+6$ and $k$ is odd} \\
                      \end{array}
               \right.
 $$
and
$$
 |D_n|=            \left\{\begin{array}{ll}
                             2+14\cdot k/2+2\cdot 2            & \mbox{\hspace{0.2cm} if $n=4k+3$ and $k$ is even} \\
                             2+14\cdot k/2+3\cdot 2            & \mbox{\hspace{0.2cm} if $n=4k+4$ and $k$ is even} \\
                             2+14\cdot k/2+4\cdot 2            & \mbox{\hspace{0.2cm} if $n=4k+5$ and $k$ is even} \\
                             \vspace{10pt}
                             2+14\cdot k/2+5\cdot 2            & \mbox{\hspace{0.2cm} if $n=4k+6$ and $k$ is even} \\
                             2+14\cdot (k-1)/2+7+2\cdot 2        & \mbox{\hspace{0.2cm} if $n=4k+3$ and $k$ is odd} \\
                             2+14\cdot (k-1)/2+7+3\cdot 2        & \mbox{\hspace{0.2cm} if $n=4k+4$ and $k$ is odd} \\
                             2+14\cdot (k-1)/2+7+4\cdot 2        & \mbox{\hspace{0.2cm} if $n=4k+5$ and $k$ is odd} \\
                             2+14\cdot (k-1)/2+7+5\cdot 2        & \mbox{\hspace{0.2cm} if $n=4k+6$ and $k$ is odd} \\
                      \end{array}
               \right.
     =\lceil\frac{7n+3}{4}\rceil.
 $$
Hence we have $\g_2(G_{4,n})\leq |D_n|=\lceil\frac{7n+3}{4}\rceil.$

\begin{thm}\label{t4}
For any positive integer $n\geq 3$,
 $
 \g_2(G_{4,n})=\lceil\frac{7n+3}{4}\rceil.
 $
\end{thm}

\begin{pf}  For $n\in \{3,4,5,6\}$, we can check straight that the
result is true. In the following, assume that $n\geq 7$. We merely prove that
$\g_2(G_{4,n})\geq \lceil\frac{7n+3}{4}\rceil$.

By Lemma \ref{l4}, $G_{4,n}$ must contain at least one
$\g_2(G_{4,n})$-set with $N_1=0$. So we can choose a
$\g_2(G_{4,n})$-set $D$ from $V(G_{4,n})$ such that $N_0(D)=0$ and
$N_2(D)$ is as large as possible. Then we have
\begin{equation}\label{e4.5}
    N_1(D)+N_2(D)+N_3(D)+N_4(D)=n
\end{equation}
and
\begin{equation}\label{e4.6}
    N_1(D)+2N_2(D)+3N_3(D)+4N_4(D)=|D|=\g_2(G_{4,n})\leq \lceil\frac{7n+3}{4}\rceil.
\end{equation}
\begin{claim}\label{c1}
 $c_1(D)=c_2(D)=c_{n-1}(D)=c_n(D)=2.$
\end{claim}

\noindent \textbf{Proof of Claim \ref{c1}.}\ \ By the symmetry of
$G_{4,n}$, we merely prove $c_1(D)=c_2(D)=2$. Supposed that
$c_1(D)\neq 2$ or $c_2(D)\neq 2$, then we have either $c_1(D)=3$ and
$c_2(D)=1$ or $c_1(D)+c_2(D)\geq 5$ by Lemma \ref{l2} and the
definition of $D$.

If $c_1(D)=3$ and $c_2(D)=1$, then we can replace one vertex of
$\mathscr{C}_1(D)$ by one of $C_2\setminus \mathscr{C}_2(D)$ and get
a new $\g_2(G_{4,n})$-set, denoted by $S_1$, which satisfies
$N_0(S_1)=0$ and $N_2(S_1)\geq N_2(D)+2$. This is a contradiction
with the choice of $D$.

If $c_1(D)+c_2(D)\geq 5$ and $c_3(D)=1$, then there exists a
$\g_2(G_{4,n})$-set $S_2$ such that $c_1(S_2)=c_2(S_2)=c_3(S_2)=2$,
$ \mathscr{C}_3(S_2)\supseteq \mathscr{C}_3(D)$ and
$\mathscr{C}_j(S_2)=\mathscr{C}_j(D)$ ($4\leq j\leq n$). Since
$N_0(S_2)=0$ and $N_2(S_2)>N_2(D)$, we get a contradiction with the
choice of $D$.

If $c_1(D)+c_2(D)\geq 5$ and $c_3(D)\geq 2$, then we can obtain a
new $2$-dominating set $S_3$ of $G_{4,n}$ from $D$ by replacing
$\mathscr{C}_1(D)$ and $\mathscr{C}_2(D)$ by two vertices in $C_1$
and $C_2$, respectively. Obviously, that $|S_3|<|D|=\g_2(G_{4,n})$
is a contradiction. This completes the proof of the claim.

Delete all $2$'s from the sequence $(c_1(D),c_2(D),\cdots,c_n(D))$,
we can obtain some subsequences, denoted by
$\mathscr{S}_1,\cdots,\mathscr{S}_h$, each of which consists of $1$,
$3$ or $4$. Notice that there exists at least one $2$'s between two
adjacent subsequences. By Claim \ref{c1}, we have
\begin{equation}\label{e4.6}
    N_2(D)\geq 2+(h-1)+2=h+3.
\end{equation}
For any $i\in [h]$, we use $\ell^1_i$, $\ell^2_i$ and $\ell^3_i$ to
denote the numbers of $1$'s, $3$'s and $4$'s in $\mathscr{S}_i$,
respectively. By Lemma \ref{l3}, we must have $\ell^3_i+\ell^4_i\geq
1$ and $\ell^1_i\leq \ell^3_i+\ell^4_i+1$, and so
\begin{equation}\label{e4.7}
    N_1(D)=\Sigma_{i=1}^h\ell^1_i\leq \Sigma_{i=1}^h(\ell^3_i+\ell^4_i+1)=N_3(D)+N_4(D)+h
\end{equation}
and
\begin{equation}\label{e4.8}
    N_3(D)+N_4(D)=\Sigma_{i=1}^h(\ell^3_i+\ell^4_i)\geq h.
\end{equation}
From (\ref{e4.5})-(\ref{e4.8}), we know that
 $$\begin{array}{rl}
 \lceil\frac{7n+3}{4}\rceil & \geq N_1(D)+2N_2(D)+3N_3(D)+4N_4(D)\\
                            & =n+N_2(D)+2(N_3(D)+N_4(D))+N_4(D)\\
                            & \geq n+(h+3)+2h+N_4(D)\\
                            & =3h+n+3+N_4(D),
 \end{array}$$
which implies $h\leq \frac{1}{3}\lceil\frac{3n+3}{4}\rceil-1-\frac{1}{3}N_4(D)$, and so
$$\begin{array}{rl}
 \lceil\frac{7n+3}{4}\rceil  \geq \g_2(G_{4,n})&= N_1(D)+2N_2(D)+3N_3(D)+4N_4(D)\\
                            & =2n+(N_3(D)+ N_4(D)-N_1(D))+N_4(D)\\
                            & \geq 2n-h+N_4(D)\\
                            & \geq 2n-\frac{1}{3}\lceil\frac{3n+3}{4}\rceil+1+\frac{4}{3}N_4(D)\\
                            & = \lceil\frac{7n+3}{4}\rceil+\frac{4}{3}N_4(D)-\left\{\begin{array}{ll}
                             1/3                       & \mbox{\hspace{0.2cm} if $n\equiv 0$ (mod $4$)} \\
                            2/3                      & \mbox{\hspace{0.2cm} if $n\equiv 1$ (mod $4$)} \\
                            1                       & \mbox{\hspace{0.2cm} if $n\equiv 2$ (mod $4$)} \\
                            0                       & \mbox{\hspace{0.2cm} if $n\equiv 3$ (mod $4$)} \\
                      \end{array}
               \right..
 \end{array}$$
This forces $N_4(D)=0$ and $\lceil\frac{7n+3}{4}\rceil-1\leq
\g_2(G_{4,n})\leq \lceil\frac{7n+3}{4}\rceil$. Supposed that
$\g_2(G_{4,n})=\lceil\frac{7n+3}{4}\rceil-1$, then $n\equiv 2$ (mod
$4$), $N_3(D)-N_1(D)=-h$ and
$h=\frac{1}{3}\lceil\frac{3n+3}{4}\rceil-1$ (which forces
$N_2(D)=h+3$ and $N_3(D)=h$), Hence
 $$
 n=N_1(D)+N_2(D)+N_3(D)=2h+(h+3)+h=4h+3,
 $$
which is a contradiction with $n\equiv 2$ (mod $4$). Thus
$\g_2(G_{4,n})=\lceil\frac{7n+3}{4}\rceil$.\end{pf}

\begin{thm}\label{t5}
For any positive integer $n\geq 7$,
 $$
 b_2(G_{4,n})=\left\{\begin{array}{ll}
            1                          & \mbox{{\rm if} {\rm $n\equiv 3$ (mod $4$)}}\\
            2                          & {\rm otherwise.}
        \end{array} \right.
  $$
\end{thm}

\begin{pf}  Since
$n\geq 7$, there are two integers $k\geq 1$ and $r\in \{3,4,5,6\}$
such that $n=4k+r$. We consider the following two cases.

\noindent \textbf{Case 1.}\ \ $n=4k+3$.

Let $e=(2,3)(3,3)$ and
$G=G_{4,n}-e$. We will prove $\g_2(G)> \g_2(G_{4,n})$, and so $b_2(G_{4,n})=1$.

Let $D$ be a $\g_2(G)$-set such that $N_0(D)=0$ and $N_2(D)$ is as
large as possible. An argument similar to that described in Claim
\ref{c1} of Theorem \ref{t4} shows $c_1(D)=c_2(D)=2$. Denote $G-C_1$
by $G'$, then $G'\cong G_{4,4k+2}-(2,2)(3,2)$.

\noindent \textbf{Subcase 1.1}\ \ $D\cap V(G')$ is a $2$-dominating
set of $G'$. By Theorem \ref{t4}, we must have
 $$
 \begin{array}{rl}
 \g_2(G)&=c_1(D)+|D\cap V(G')|\\
  & \geq 2+\g_2(G')\\
  & \geq 2+\g_2(G_{4,4k+2})=2+\lceil\frac{7(4k+2)+3}{4}\rceil=7k+7\\
  & > 7k+6=\g_2(G_{4,4k+3})\\
   &=\g_2(G_{4,n}).
 \end{array}
 $$

\noindent \textbf{Subcase 1.2}\ \ $D\cap V(G')$ isn't a
$2$-dominating set of $G'$.

If $|D\cap \{(1,3),(2,3)\}|\geq 1$ and $|D\cap \{(3,3),(4,3)\}|\geq
1$, then there exist $\mathscr{C}_1'\subseteq C_1$ and
$\mathscr{C}_2'\subseteq C_2$ with
$|\mathscr{C}_1'|=|\mathscr{C}_1'|=2$ such that
 $$
 D'=\mathscr{C}_1'\cup \mathscr{C}_2'\cup\mathscr{C}_3(D)\cup \cdots\cup\mathscr{C}_n(D)
 $$
is a new $\g_2(G)$-set and $D'\cap V(G')$ is a $2$-dominating set of
$G'$. By \textbf{Subcase 1.1}, $\g_2(G)=|D'|>\g_2(G_{4,n})$.

If $|D\cap \{(1,3),(2,3)\}|=0$ or $|D\cap \{(3,3),(4,3)\}|=0$,
without loss of generality, say $|D\cap \{(3,3),(4,3)\}|=0$, then,
to $2$-dominate $C_1$, $C_2$ and $C_3$, $\mathscr{C}_1(D)\cup
\mathscr{C}_2(D)\cup \mathscr{C}_3(D)$ must be the set of black
vertices in Figure \ref{f4}.

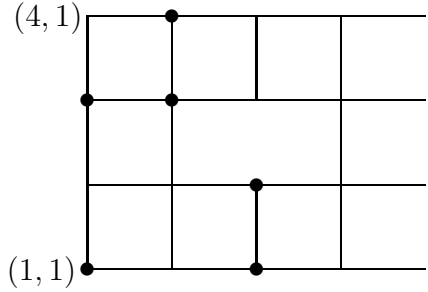
\begin{figure}[h]
\begin{center}
\psset{unit=0.8pt}

\begin{pspicture}(400,200)(0,50)

\put(120,70){\line(1,0){160}}\put(120,110){\line(1,0){160}}\put(120,150){\line(1,0){160}}\put(120,190){\line(1,0){160}}

\put(120,70){\line(0,1){120}}\put(160,70){\line(0,1){120}}\put(200,70){\line(0,1){40}}\put(200,150){\line(0,1){40}}\put(240,70){\line(0,1){120}}

\put(120,70){\circle*{6}}\put(120,150){\circle*{6}}
\put(160,150){\circle*{6}}\put(160,190){\circle*{6}}
\put(200,70){\circle*{6}}\put(200,110){\circle*{6}}
\put(82,65){$(1,1)$}\put(85,185){$(4,1)$}

 \end{pspicture}
\caption{\label{f4} \footnotesize
$\mathscr{C}_1(D)\cup\mathscr{C}_1(D)\cup\mathscr{C}_1(D)$ is the
set of black vertices}
\end{center}
\end{figure}

An argument similar to that described in the proof of Theorem
\ref{t4} shows that $N_i(D)$ ($1\leq i\leq 4$) meets (\ref{e4.5}),
(\ref{e4.7}), (\ref{e4.8}) and $N_2(D)\geq 3+(h-1)+2=h+4$ ($h$ is
defined in the proof of Theorem \ref{t4}). Hence
$$
 \begin{array}{rl}
 \g_2(G)=\Sigma_{i=1}^4iN_i(D)&=n+N_2(D)+2(N_3(D)+N_4(D))+N_4(D)\\
  & \geq n+(h+4)+2h+N_4(D)\\
  & =3h+n+4+N_4(D),
 \end{array}
 $$
and so
$$
 \begin{array}{rl}
 \g_2(G)=\Sigma_{i=1}^4iN_i(D)&=2n+(N_3(D)+N_4(D)-N_1(D))+N_4(D)\\
  & \geq 2n-h+N_4(D)\\
  & \geq 2n-\frac{1}{3}(\g_2(G)-n-4-N_4(D))+N_4(D)\\
  & =\frac{7n+4}{3}+\frac{4}{3}N_4(D)-\frac{1}{3}\g_2(G).\\
  \end{array}
 $$
This forces
$$
 \begin{array}{rl}
 \g_2(G)\geq \lceil\frac{7n+4}{4}+N_4(D)\rceil=7k+7+N_4(D)>7k+6=\g_2(G_{4,n}).\\
  \end{array}
 $$

\noindent \textbf{Case 2.}\ \ $n=4k+r$, $r\in \{4,5,6\}$.

We first show that $b_2(G_{4,n})> 1$, that is, we will prove
$\g_2(G_{4,n}-e)=\g_2(G_{4,n})$ for any edge $e\in E(G_{4,n})$.  Let
$e=(i,j)(i',j')$ and $H=G_{4,n}-e$. By the symmetry of $G_{4,n}$, we
assume that $i=1,2$ and $1\leq j\leq \lfloor\frac{n}{2}\rfloor$.

Recalled the definitions of $T_{4k+3}$ and $D_{4k+3}$, we know that
the set $D_{4k+3}$ of black vertices in block
 $$T_{4k+3}=\left\{\begin{array}{ll}
                 B_1(B_3B_4)^{\frac{k}{2}}B_2B_1
                 & \mbox{\hspace{0.2cm} if $k$ is even}\\
                 B_1(B_3B_4)^{\frac{k-1}{2}}B_3B_2B_1
                 & \mbox{\hspace{0.2cm} if $k$ is odd}
                 \end{array}
            \right.
                 $$
is a $\g_2(G_{4,4k+3})$-set by Theorem \ref{t4}, and we can
construct an edge subset $E_{4k+3}$ of $G_{4,4k+3}$ as follows:
 $$
E_{4k+3}=\{xy\in E(G_{4,4k+3}) \ |\ x\in D_{4k+3}, y\notin D_{4k+3} \mbox{
and } |N_{G_{4,4k+3}}(y)\cap D_{4k+3}|=2\}.
 $$
Illustrations of the block $T_{4k+3}$, the $\g_2(G_{4k+3})$-set
$D_{4k+3}$ and the constructed edge subset $E_{4k+3}$ are shown in
Figure \ref{f5}.

\begin{figure}[h]
\begin{center}
\psset{unit=0.8pt}

\begin{pspicture}(400,270)(0,80)

\put(-38,260){\line(1,0){218}} \put(-40,230){\line(1,0){220}}
\put(-40,200){\line(1,0){220}} \put(-38,170){\line(1,0){218}}

\put(235,260){\line(1,0){198}} \put(235,230){\line(1,0){200}}
\put(235,200){\line(1,0){200}} \put(235,170){\line(1,0){198}}

\put(200,215){$\cdots$}
 \put(-55,150){$(1,1)$}
\put(-55,270){$(4,1)$}
\put(405,150){$(1,4k+3)$}\put(405,270){$(4,4k+3)$}

\put(-40,172){\line(0,1){86}} \put(-20,170){\line(0,1){90}}
\put(0,170){\line(0,1){90}} \put(20,170){\line(0,1){90}}
\put(40,170){\line(0,1){90}} \put(60,170){\line(0,1){90}}
\put(80,170){\line(0,1){90}}\put(100,170){\line(0,1){90}}
\put(120,170){\line(0,1){90}}\put(140,170){\line(0,1){90}}

\put(435,172){\line(0,1){86}}\put(415,170){\line(0,1){90}}
\put(395,170){\line(0,1){90}}\put(375,170){\line(0,1){90}}
\put(355,170){\line(0,1){90}}\put(335,170){\line(0,1){90}}
\put(315,170){\line(0,1){90}}\put(295,170){\line(0,1){90}}
\put(275,170){\line(0,1){90}}

\linethickness{0.5mm}
\put(-40,172){\line(0,1){28}}\put(-40,230){\line(0,1){28}}
\put(0,200){\line(0,1){30}}\put(40,200){\line(0,1){30}}\put(80,200){\line(0,1){30}}\put(120,200){\line(0,1){30}}
\put(275,200){\line(0,1){30}}\put(315,200){\line(0,1){30}}\put(355,200){\line(0,1){30}}\put(395,200){\line(0,1){30}}
\put(435,172){\line(0,1){28}}\put(435,230){\line(0,1){28}}

\put(60,170){\line(0,1){30}}\put(140,170){\line(0,1){30}}\put(335,170){\line(0,1){30}}\put(415,170){\line(0,1){30}}

\put(60,230){\line(0,1){30}}\put(140,230){\line(0,1){30}}\put(335,230){\line(0,1){30}}\put(-20,230){\line(0,1){30}}

\put(-38,170){\line(1,0){18}}\put(-38,260){\line(1,0){18}}\put(335,170){\line(1,0){98}}\put(415,260){\line(1,0){18}}
\put(-40,230){\line(1,0){20}}\put(0,230){\line(1,0){20}}\put(20,230){\line(1,0){20}}\put(60,230){\line(1,0){20}}\put(120,230){\line(1,0){20}}
\put(60,170){\line(1,0){80}}
\put(-20,260){\line(1,0){80}}\put(295,260){\line(1,0){40}}
\put(40,200){\line(1,0){20}}\put(80,200){\line(1,0){40}}
\put(275,230){\line(1,0){40}}\put(335,230){\line(1,0){20}}
\put(315,200){\line(1,0){20}}\put(355,200){\line(1,0){40}}\put(415,200){\line(1,0){20}}

\put(-40,230){\circle*{6}} \put(-40,200){\circle*{6}}
\put(-20,260){\circle*{6}} \put(-20,170){\circle*{6}}
\put(0,200){\circle*{6}} \put(20,230){\circle*{6}}
\put(20,260){\circle*{6}} \put(20,170){\circle*{6}}
\put(40,200){\circle*{6}} \put(60,170){\circle*{6}}
\put(60,260){\circle*{6}} \put(80,230){\circle*{6}}
\put(100,170){\circle*{6}}\put(100,200){\circle*{6}}
\put(100,260){\circle*{6}}\put(120,230){\circle*{6}}
\put(140,170){\circle*{6}}\put(140,260){\circle*{6}}

\put(275,200){\circle*{6}}\put(295,170){\circle*{6}}
\put(295,230){\circle*{6}}
\put(295,260){\circle*{6}}\put(315,200){\circle*{6}}
\put(335,170){\circle*{6}}
\put(335,260){\circle*{6}}\put(355,230){\circle*{6}}
\put(375,170){\circle*{6}}\put(375,200){\circle*{6}}
\put(375,260){\circle*{6}}\put(395,230){\circle*{6}}
\put(415,170){\circle*{6}}\put(415,260){\circle*{6}}
\put(435,200){\circle*{6}}\put(435,230){\circle*{6}}

\put(-40,170){\circle{4}}\put(-40,260){\circle{4}}
\put(435,170){\circle{4}}\put(435,260){\circle{4}}

\put(160,120){\scriptsize $(a)$\ \ When $k$ is even}
\end{pspicture}

\begin{pspicture}(400,270)(0,80)

\put(-38,260){\line(1,0){198}} \put(-40,230){\line(1,0){200}}
\put(-40,200){\line(1,0){200}} \put(-38,170){\line(1,0){198}}

\put(195,260){\line(1,0){258}} \put(195,230){\line(1,0){260}}
\put(195,200){\line(1,0){260}} \put(195,170){\line(1,0){258}}

\put(170,215){$\cdots$}
 \put(-55,150){$(1,1)$}
\put(-55,270){$(4,1)$}
\put(410,150){$(1,4k+3)$}\put(410,270){$(4,4k+3)$}

\put(-40,172){\line(0,1){86}} \put(-20,170){\line(0,1){90}}
\put(0,170){\line(0,1){90}} \put(20,170){\line(0,1){90}}
\put(40,170){\line(0,1){90}} \put(60,170){\line(0,1){90}}
\put(80,170){\line(0,1){90}}\put(100,170){\line(0,1){90}}
\put(120,170){\line(0,1){90}}\put(140,170){\line(0,1){90}}

\put(455,172){\line(0,1){86}}
\put(435,170){\line(0,1){90}}\put(415,170){\line(0,1){90}}
\put(395,170){\line(0,1){90}}\put(375,170){\line(0,1){90}}
\put(355,170){\line(0,1){90}}\put(335,170){\line(0,1){90}}
\put(315,170){\line(0,1){90}}\put(295,170){\line(0,1){90}}
\put(275,170){\line(0,1){90}}\put(255,170){\line(0,1){90}}
\put(235,170){\line(0,1){90}} \put(215,170){\line(0,1){90}}

\linethickness{0.5mm}
\put(-40,172){\line(0,1){28}}\put(-40,230){\line(0,1){28}}
\put(0,200){\line(0,1){30}}\put(40,200){\line(0,1){30}}\put(80,200){\line(0,1){30}}\put(120,200){\line(0,1){30}}
\put(215,200){\line(0,1){30}}\put(255,200){\line(0,1){30}}\put(295,200){\line(0,1){30}}\put(335,200){\line(0,1){30}}
\put(375,200){\line(0,1){30}}\put(415,200){\line(0,1){30}}\put(455,172){\line(0,1){28}}\put(455,230){\line(0,1){28}}

\put(-38,170){\line(1,0){18}}\put(-38,260){\line(1,0){18}}\put(275,170){\line(1,0){80}}\put(415,260){\line(1,0){18}}\put(435,170){\line(1,0){18}}
\put(-40,230){\line(1,0){20}}\put(0,230){\line(1,0){20}}\put(20,230){\line(1,0){20}}\put(60,230){\line(1,0){20}}\put(120,230){\line(1,0){20}}
\put(60,170){\line(1,0){80}}
\put(-20,260){\line(1,0){80}}\put(215,260){\line(1,0){60}}\put(355,260){\line(1,0){98}}
\put(40,200){\line(1,0){20}}\put(80,200){\line(1,0){40}}
\put(215,230){\line(1,0){40}}\put(275,230){\line(1,0){20}}\put(335,230){\line(1,0){20}}\put(375,230){\line(1,0){40}}\put(435,230){\line(1,0){20}}
\put(255,200){\line(1,0){20}}\put(295,200){\line(1,0){40}}\put(355,200){\line(1,0){20}}

\put(60,170){\line(0,1){30}}\put(140,170){\line(0,1){30}}\put(275,170){\line(0,1){30}}\put(355,170){\line(0,1){30}}

\put(60,230){\line(0,1){30}}\put(140,230){\line(0,1){30}}\put(275,230){\line(0,1){30}}\put(-20,230){\line(0,1){30}}\put(355,230){\line(0,1){30}}
\put(435,230){\line(0,1){30}}

\put(-40,230){\circle*{6}} \put(-40,200){\circle*{6}}
\put(-20,260){\circle*{6}} \put(-20,170){\circle*{6}}
\put(0,200){\circle*{6}} \put(20,230){\circle*{6}}
\put(20,260){\circle*{6}} \put(20,170){\circle*{6}}
\put(40,200){\circle*{6}} \put(60,170){\circle*{6}}
\put(60,260){\circle*{6}} \put(80,230){\circle*{6}}
\put(100,170){\circle*{6}}\put(100,200){\circle*{6}}
\put(100,260){\circle*{6}}\put(120,230){\circle*{6}}
\put(140,170){\circle*{6}}\put(140,260){\circle*{6}}

\put(215,200){\circle*{6}}
\put(235,170){\circle*{6}}\put(235,230){\circle*{6}}
\put(235,260){\circle*{6}} \put(255,200){\circle*{6}}
\put(275,170){\circle*{6}}\put(275,260){\circle*{6}}
\put(295,230){\circle*{6}} \put(315,200){\circle*{6}}
\put(315,260){\circle*{6}} \put(315,200){\circle*{6}}
\put(315,170){\circle*{6}} \put(335,230){\circle*{6}}
\put(355,170){\circle*{6}} \put(355,260){\circle*{6}}
\put(375,200){\circle*{6}}
\put(395,170){\circle*{6}}\put(395,230){\circle*{6}}
\put(395,260){\circle*{6}} \put(415,200){\circle*{6}}
\put(435,170){\circle*{6}}\put(435,260){\circle*{6}}
\put(455,200){\circle*{6}}\put(455,230){\circle*{6}}

\put(-40,170){\circle{4}}\put(-40,260){\circle{4}}
\put(455,170){\circle{4}}\put(455,260){\circle{4}}

 \put(160,120){\scriptsize $(b)$\ \ When $k$ is odd}
\end{pspicture}

\caption{\label{f5} \footnotesize Illustrations of the block $T_{4k+3}$,
the subset $D_{4k+3}$ composed of all black vertices and the edge subset $E_{4k+3}$ consisting of all black edges}
\end{center}
\end{figure}

Obviously, for any $e\notin E_{4k+3}$, $D_{4k+3}$ is still a
$2$-dominating set of $G_{4,4k+3}-e$. Note that $G_{4,4k+3}$ is an
induced subgraph of $G_{4,n}$. So $E_{4k+3}\subseteq E(G_{4,n})$.
By the definitions of $D_{4k+3}$ and $D_n$, $D_{4k+3}\subseteq D_n$,
and so $D_{n}$ is also a $2$-dominating set of $H=G_{4,n}-e$. This
implies $\g_2(H)=\g_2(G_{4,n})$. In the following, let $e\in
E_{4k+3}$. As assuming that $e=(i,j)(i',j')$ satisfies $i=1,2$ and
$1\leq j\leq \lfloor\frac{n}{2}\rfloor$, we merely consider four
cases.

If $e=(1,1)(2,1)$, then it is easy to see that
 $$
 D_n'=(D_n\setminus \{(2,1),(1,2)\})\cup \{(1,1),(2,2)\}
 $$
 is a $2$-dominating set of
$H$ and $\g_2(H)\leq |D_n'|=|D_n|=\g_2(G_{4,n})$, which implies $\g_2(H)=\g_2(G_{4,n})$.

If $e=(1,1)(1,2)$, then it is obvious that
 $$
 D_n'=(D_n\setminus \{(2,1)\})\cup \{(1,1)\}
 $$
 is a $2$-dominate
of $H$ with $|D_n'|=|D_n|=\g_2(G_{4,n})$. So $\g_2(H)=\g_2(G_{4,n})$.

If $e=(i,j)(i+1,j)$ and $j\geq 2$, then obviously edge
$(i,j-1)(i+1,j-1)\notin E_{4k+3}$ by Figure \ref{f5}.
Hence for $n=4k+r$ $(4\leq r\leq 6)$, the set of black vertices in the block
 $$
T_n'= \left\{\begin{array}{ll}
            B_2T_{4k+3}                                   & \mbox{if $r=4$}\\
            B_2T_{4k+3}B_2                 & \mbox{if $r=5$}\\
            B_2T_{4k+3}B_2B_1                 & \mbox{if $r=6$}
             \end{array}
       \right.
 $$
 is a $2$-dominating set of $H$ with  cardinality
$$\g_2(G_{4,4k+3})+\left\{\begin{array}{ll}
 2   & \mbox{if $r=4$}\\
 4    & \mbox{if $r=5$}\\
 6    & \mbox{if  $r=6$}
 \end{array}
 \right.
 =\lceil\frac{7n+3}{4}\rceil=\g_2(G_{4,n}).$$

If $e=(i,j)(i,j+1)$ and $j\geq 2$, then we can find edge
$(5-i,j)(5-i,j+1)\notin E_{4k+3}$ by observing Figure \ref{f5}. For
$n=4k+r$ and $r\in \{4,5,6\}$, since the black vertex set $D_n$ in
the block $T_n$ is a $\g_2(G_{4,n})$-set by Theorem \ref{t4}, set
 $$
D_n'=\{(i,j)\ |\ (5-i,j)\in D_n\}
 $$
$2$-dominate $H$, and so $\g_2(H)=|D_n'|=|D_n|=\g_2(G_{4,n})$.

To the end, we only need to prove that $b_2(G_{4,n})\leq 2$. Let
 $$e_1=(1,5)(1,6),\ e_2=(4,5)(4,6) \mbox{\ \ and\ \ } L=G_{4,n}-\{e_1,e_2\}.$$
 We will
show $\g_2(L)> \g_2(G_{4,n})$, and so $b_2(G_{4,n})\leq 2$. Denote the subgraph
$L-\cup_{j=1}^5C_j$ by $L_1$.

Let $S$ be a $\g_2(L)$-set such that $N_0(S)=0$ and $N_2(S)$ is as
large as possible. An argument similar to that described in Claim
\ref{c1} of Theorem \ref{t4} shows $c_1(S)=c_2(S)=2$. Since
 $
 deg_L((1,5))=deg_L((4,5))=2,
 $
 to $2$-dominate $(1,5)$ and $(4,5)$, we have
  $$
  \mbox{$|S\cap \{(1,5),(2,5)\}|\geq 1$ \mbox{ and } $|S\cap \{(3,5),(4,5)\}|\geq
  1$,}
  $$
which implies $c_5(S)\geq 2$. By Lemma \ref{l3}, if $c_3(S)=1$
(resp. $c_4(S)=1$) then $c_2(S)+c_3(S)+c_4(S)\geq 6$ (resp.
$c_3(S)+c_4(S)+c_5(S)\geq 6$). Hence we have
$$c_1(S)+c_2(S)+c_3(S)+c_4(S)+c_5(S)\geq 10.$$

If $S\cap V(L_1)$ is a $2$-dominating set of $L_1$, then by $L_1\cong G_{4,n-5}$, we have
 $$\begin{array}{ll}
 \g_2(L)=|S|&=\sum_{j=1}^5c_j(S)+|S\cap V(L_1)|\\
 &\geq 10+\lceil\frac{7(n-5)+3}{4}\rceil\\
 &> \lceil\frac{7n+3}{4}\rceil\\
 &=\g_2(G_{4,n}).
 \end{array}
 $$

If $S\cap V(L_1)$ isn't a $2$-dominating set of $L_1$, then at least
one of $(2,6)$ and $(3,6)$ can't be $2$-dominated by $S\cap V(L_1)$.
Since
$$
  \mbox{$|S\cap \{(1,6),(2,6)\}|\geq 1$ \mbox{ and } $|S\cap \{(3,6),(4,6)\}|\geq
  1$}
  $$
by $deg_L((1,6))=deg_L((4,6))=2$, we must have $$S\cap
\{(1,6),(2,6)\}=\{(1,6)\}\ \mbox{and } S\cap
\{(3,6),(4,6)\}=\{(4,6)\}.$$ Thus set
 $$\begin{array}{rr}
 \{(2,1),(3,1)\}\cup\{(1,2),(4,2)\}\cup\{(2,3),(3,3)\}\cup\{(1,4),(4,4)\}\\
 \cup\{(2,5),(3,5)\}\cup (S\cap V(L_1))
 \end{array}$$
is also a $2$-dominating set of $L$ and set
 $$
 \{(2,5),(3,5)\}\cup (S\cap V(L_1))
 $$
is a $2$-dominating set of subgraph $G_{4,n}-\cup_{j=1}^4C_j$. Note
that $G_{4,n}-\cup_{j=1}^4C_j$ is isomorphic to $G_{4,n-4}$. Hence
by Theorem \ref{t4}, we must have
 $$
 \begin{array}{rl}
\g_2(L)=|S|   &=\sum_{j=1}^5c_j(S)+|S\cap V(L_1)|\\
              &\geq 10+|S\cap V(L_1)|\\
              &=8+|\{(2,5),(3,5)\}\cup (S\cap V(L_1))|\\
              &\geq 8+\lceil\frac{7(n-4)+3}{4}\rceil\\
              &=\lceil\frac{7n+3}{4}\rceil+1\\
              &>\g_2(G_{4,n}).
\end{array}
 $$
The proof of the theorem is completed.\end{pf}


\begin{thebibliography}{99}


\bibitem{bcf05}
M. Blidia and M. Chellali, O. Favaron, Independence and 2-domination
in trees. {\it Austral. J. Combin.}, 33(2005), 317-327.


\bibitem{bcv06}
M. Blidia, M. Chellali and L. Volkmann, Some bounds on the
$p$-domination number in trees. {\it Discrete Math.}, 306(2006),
2031-2037.

\bibitem{bcv06b}
M. Blidia, M. Chellali, and L. Volkmann, Bounds of the 2-domination
number of graphs. {\it Util. Math}. 71 (2006), 209-216.



\bibitem{cd06}
K. Carlson and M. Develin, On the bondage number of planar and
directed graphs. {\it Discrete Math.}, {\bf 306} (8-9) (2006),
820-826.

\bibitem{cr90}
Y. Caro and Y. Roditty, A note on the $k$-domination number of a
graph. {\it Internat. J. Math. Sci.}, 13(1990), 205-206.


\bibitem{cc93}
T. Y. Chang, and W. E. Clark, The domination number of the $5\times
n$ and $6\times n$ grids graphs. {\it J. Graph Theory}, 17(1993),
81-108.


\bibitem{cch94}
T. Y. Chang, W. E. Clark and E. O. Hare, Domination numbers of
complete grid graphs, I. {\it Ars Combin.}, 38 (1994), 97-111.



\bibitem{c06}
M. Chellali, Bounds on the 2-domination number in cactus graphs.
{\it Opuscula Math}. 26 (1)(2006), 5-12.


\bibitem{dlpw10}
E. DeLaVi\~na, C. E. Larson, R. Pepper, B. Waller, Graffiti.pc on
the 2-domination number of a graph. Proceedings of the Forty-First
Southeastern International Conference on Combinatorics, Graph Theory
and Computing. {\it Congressus Numerantium}, 203 (2010), 15-32.

\bibitem{dghpv11}
E. DeLaVi\~na, W. Goddard, M. A. Henning, R. Pepper, E. R. Vaughan,
Bounds on the $k$-domination number of a graph. {\it Applied
Mathematics Letters}, 24 (6) (2011), 996-998.

\bibitem{f85}
O. Favaron, On a conjecture of Fink and Jacobson concerning
$k$-domination and $k$-dependence. {\it J. Combin. Theory, Ser. B},
39(1985), 101-102.


\bibitem{fj851}
J. F. Fink and M. S. Jacobson, $n$-domination in graphs. In {\it
Graph Theory with Applications to Algorithms and Computer Science}
(Y. Alavi, A. J. Schwenk eds), 283-300, Wiley, New York, (1985).

\bibitem{fj852}
J. F. Fink and M. S. Jacobson, On $n$-domination, $n$-dependence and
forbidden subgraphs. Graph Theory with Applications to Algorithms
and Computer Science, 301-311, John Wiley and Sons, New York,
(1985).


\bibitem{fjkr90}
J. F. Fink, M. S. Jacobson, L. F. Kinch and J. Roberts, The bondage
number of a graph. {\it Discrete Math.}, 86(1990), 47-57.


\bibitem{frv03}
M. Fischermann, D. Rautenbach and L. Volkmann, Remarks on the
bondage number of planar graphs. {\it Discrete Mathematics}, {\bf
260} (2003), 57-67.


\bibitem{gprt11}
D. Gon\c{c}alves, A. Pinlon, M. Rao and S. Thomasse, The domination
number of grids, 
arXiv:1102.5206v1[cs.DM], 25 Feb. 2011.


\bibitem{hv07}
A. Hansberg, L. Volkmann, Characterization of block graphs with
equal 2-domination number and domination number plus one. {\it
Discuss. Math. Graph Theory}, 27 (1)(2007), 93-103.


\bibitem{hv08a}
A. Hansberg, L. Volkmann, Characterization of unicyclic graphs with
equal 2-domination number and domination number plus one. {\it
Utilitas Mathematica},
77 (2008), 265-276.

\bibitem{hv08b}
A. Hansberg, L. Volkmann, On graphs with equal domination and
2-domination numbers. {\it Discrete Math}. 308 (11)(2008),
2277-2281.


\bibitem{hr92}
B. L. Hartnell and D. F. Rall, A characterization of trees in which
no edge is essential to the domination number. {\it Ars Comb.},
33(1992), 65-76.

\bibitem{hl94}
B. L. Hartnell and D. F. Rall, Bounds on the bondage number of a
graph. {\it Discrete Math.}, 128(1994), 173-177.

\bibitem{hs981}
T. W. Haynes, S. T. Hedetniemi and P. J. Slater, {\it Fundamentals
of Domination in Graphs}. New York, Marcel Deliker, (1998).

\bibitem{hs982}
T. W. Haynes, S. T. Hedetniemi and P. J. Slater, {\it Domination in
Graphs: Advanced Topics}. New York, Marcel Deliker (1998).


\bibitem{hcx09}
F.-T. Hu, Y.-C. Cao and J.-M. Xu, The bondage number of mesh
networks, a manuscript submitted to {\it Comput. Math. Appl.}

\bibitem{hx06}
J. Huang and J.-M. Xu, The bondage numbers of extended de Bruijn and
Kautz digraphs. {\it Comput. Math. Appl.}, 51(6-7)(2006), 1137-1147.


\bibitem{jk83}
M. S. Jacobson and L. F. Kinch, On the domination number of product
of a graph: I. {\it Ars Comb.}, 18(1983), 33-44.


\bibitem{ksk05}
L.-Y. Kang, M. Y. Sohn and H. K. Kim, Bondage number of the discrete
torus $C_n\times C_4$. {\it Discrete Math.}, 303(2005), 80-86.


\bibitem{k10}
M. Krzywkowski, An alternative proof of a lower bound on the
2-domination number of a tree.
{\it Int. J. Mod. Math}. 5 (3)(2010), 325-326.

\bibitem{lx10}
Y. Lu and J.-M. Xu, The $p$-bondage number of trees. {\it Graphs and
Comb.}, 27(2011), 129-141.



\bibitem{s09}
R. S. Shaheen, Bounds for the 2-domination number of toroidal grid
graphs. {\it Int. J. Comput. Math}., 86 (4)(2009), 584-588.


\bibitem{t97}
U. Teschner, New results about the bondage number of a graph. {\it
Discrete Math.}, 171(1997), 249-259.


\bibitem{v08}
L. Volkmann, A Nordhaus-Gaddum-type result for the 2-domination
number. {\it J. Combin. Math. Combin. Comput.}, 64 (2008), 227-235.


\end{thebibliography}
\end{document}